\documentclass[11pt,twoside]{article}
\usepackage[utf8]{inputenc}
\usepackage{amsmath}
\usepackage{amsfonts}
\usepackage{amssymb}
\usepackage{amsthm}
\usepackage[english]{babel}
\usepackage[utf8]{inputenc}
\usepackage{hyperref}
\usepackage[alphabetic]{amsrefs}
\usepackage{geometry}
\usepackage{graphicx}
\usepackage{mathrsfs}
\usepackage{parskip}
\usepackage{paralist}

\hypersetup{
	colorlinks   = true, 
	urlcolor     = blue, 
	linkcolor    = blue,
	citecolor    = blue 
}

\usepackage{fancyhdr}
\makeatletter
\def\thm@space@setup{%
	\thm@preskip=\parskip \thm@postskip=0pt
}
\makeatother
\usepackage{etoolbox}
\AtBeginEnvironment{proof}{\vspace*{-\parskip}}

\usepackage{titlesec}
\titlespacing{\section}{0pt}{1ex}{-1ex}
\titlespacing{\subsection}{0pt}{1ex}{-1ex}

\begin{document}
\title{\textbf{Existence of multiple closed CMC hypersurfaces with small mean curvature}}
\author{Akashdeep Dey \thanks{Email: adey@math.princeton.edu, dey.akash01@gmail.com}}
\date{}
\maketitle

\theoremstyle{plain}
\newtheorem{thm}{Theorem}[section]
\newtheorem{lem}[thm]{Lemma}
\newtheorem{pro}[thm]{Proposition}
\newtheorem{clm}[thm]{Claim}
\newtheorem*{thm*}{Theorem}
\newtheorem*{lem*}{Lemma}
\newtheorem*{clm*}{Claim}

\theoremstyle{definition}
\newtheorem{defn}[thm]{Definition}
\newtheorem{ex}[thm]{Example}
\newtheorem{rmk}[thm]{Remark}

\numberwithin{equation}{section}

\newcommand{\mf}{manifold\;}\newcommand{\vf}{varifold\;}\newcommand{\mfs}{manifolds\;}\newcommand{\hy}{hypersurface\;}\newcommand{\hys}{hypersurfaces\;}\newcommand{\rim}{Riemannian manifold\;}\newcommand{\Rm}{Riemannian\;}\newcommand{\cn}{constant\;}\newcommand{\mt}{metric\;} \newcommand{\st}{such that\;}\newcommand{\cH}{\mathcal{H}}\newcommand{\Thm}{Theorem\;}\newcommand{\Eqn}{Equation\;}\newcommand{\te}{there exists\;}\newcommand{\tf}{Therefore, \;}\newcommand{\wrt}{with respect to\;}\newcommand{\bbr}{\mathbb{R}}\newcommand{\bbn}{\mathbb{N}}\newcommand{\mres}{\mathbin{\vrule height 1.6ex depth 0pt width
0.13ex\vrule height 0.13ex depth 0pt width 1.3ex}}\newcommand{\ra}{\rightarrow}\newcommand{\sps}{Suppose\;}\newcommand{\del}{\partial}\newcommand{\al}{\alpha}\newcommand{\be}{\beta}\newcommand{\ve}{\varepsilon}\newcommand{\vp}{\varphi}\newcommand{\ga}{\gamma}\newcommand{\Ga}{\Gamma}\newcommand{\de}{\delta}\newcommand{\De}{\Delta}\newcommand{\ph}{\phi}\newcommand{\Ph}{\Phi}\newcommand{\ps}{\psi}\newcommand{\Ps}{\Psi}\newcommand{\la}{\lambda}\newcommand{\La}{\Lambda}\newcommand{\si}{\sigma}\newcommand{\Si}{\Sigma}\newcommand{\om}{\omega}\newcommand{\Om}{\Omega}\newcommand{\ta}{\tau}\newcommand{\et}{\eta}\newcommand{\rh}{\rho}\newcommand{\tht}{\theta}\newcommand{\Th}{\Theta}\newcommand{\ka}{\kappa}\newcommand{\ze}{\zeta}\newcommand{\ch}{\chi}\newcommand{\seq}{sequence\;}\newcommand{\cts}{continuous\;}\newcommand{\bF}{\mathbf{F}} \newcommand{\cF}{\mathcal{F}} \newcommand{\bM}{\mathbf{M}}\newcommand{\bL}{\mathbf{L}}\newcommand{\w}{with optimal regularity\;}\newcommand{\cR}{\mathcal{R}}\newcommand{\cS}{\mathcal{S}}\newcommand{\bet}{\mathbf{\eta}}\newcommand{\scrX}{\mathscr{X}}\newcommand{\cZ}{\mathcal{Z}}

\newcommand{\norm}[1]{\left\|#1\right\|}\newcommand{\md}[1]{\left|#1\right|}\newcommand{\Md}[1]{\Big|#1\Big|}\newcommand{\db}[1]{[\![#1]\!]}
\newcommand{\vol}{\operatorname{Vol}}\newcommand{\cl}{\operatorname{Clos}}
\newcommand{\vt}{\operatorname{VarTan}}\newcommand{\spt}{\operatorname{spt}}

\newcommand{\cm}{\mathcal{C}(M)}
\newcommand{\zn}{\mathcal{Z}_n(M;\mathbf{F};\mathbb{Z}_2)}
\newcommand{\znf}{\mathcal{Z}_n(M;\mathcal{F};\mathbb{Z}_2)}
\newcommand{\ac}{\mathcal{A}^c}
\newcommand{\tx}{\tilde{X}}
\newcommand{\stx}{S\tilde{X}}
\newcommand{\tz}{\tilde{Z}}
\newcommand{\tp}{\tilde{\Phi}}
\newcommand{\ts}{\tilde{\Psi}}
\newcommand{\bbz}{\mathbb{Z}}
\newcommand{\tsia}{\tilde{\sigma}_{2i-1}}
\newcommand{\tsib}{\tilde{\sigma}_{2i}}

\vspace{-2ex}
\begin{abstract}\vspace{-1.5ex}
\noindent
Let $(M^{n+1},g)$ be a closed Riemannian manifold, $n+1\geq 3$. We will prove that for all $m \in \bbn$, \te $c^{*}(m)>0$, which depends on \(g\), \st if \(0<c<c^{*}(m)\), \((M,g)\) contains at least \(m\) many closed \(c\)-CMC hypersurfaces with optimal regularity.  More quantitatively, there exists a constant \(\ga_0\), depending on \(g\), \st for all \(c>0\), there exist at least \(\ga_0c^{-\frac{1}{n+1}}\) many closed \(c\)-CMC hypersurfaces (with optimal regularity) in \((M,g)\). This extends the theorem of Zhou and Zhu \cite{zz1}, where they proved the existence of at least one closed \(c\)-CMC hypersurface in \((M,g)\).
\end{abstract}

\section{Introduction}
Constant mean curvature (CMC) hypersurfaces are the critical points of the area functional \wrt the variations which preserve the enclosed volume. The hypersurfaces with constant mean curvature zero are called minimal hypersurfaces; they are the critical points of the area functional.

There is a well-developed existence theory for the closed minimal \hys in closed \Rm manifolds. By the combined work of Almgren \cite{alm}, Pitts \cite{pitts} and Schoen-Simon \cite{ss}, every closed \rim \((M^{n+1},g)\), \(n+1 \geq 3\), contains a closed minimal hypersurface, which is smooth and embedded outside a singular set of Hausdorff dimension \(\leq n-7\).

Since the cohomology ring of \(\mathcal{Z}_n(M^{n+1},\bbz_2)\) (with coefficients in \(\bbz_2\)) is isomorphic to the polynomial ring \(\bbz_2[\xi]\), from the finite dimensional Morse theory one expects that \((M,g)\) contains infinitely many closed minimal hypersurfaces. This was conjectured by Yau \cite{yau}. Yau's conjecture has been completely resolved when the ambient dimension \(3\leq n+1\leq 7\). By the works of Marques-Neves \cite{mn_ricci_positive} and Song \cite{song}, in every closed Riemannian manifold \((M^{n+1},g)\), \(3\leq n+1 \leq 7\), there exist infinitely many closed minimal hypersurfaces. In higher dimensions, Li \cite{li} has proved that every closed manifold \(M^{n+1}\), \(n+1 \geq 3\), equipped with a generic metric, contains infinitely many closed minimal hypersurfaces with optimal regularity.

When \(3\leq n+1 \leq 7\), stronger results have been obtained for generic set of metrics. In \cite{imn}, Irie, Marques and Neves proved that for a generic metric, the union of all closed minimal hypersurfaces is dense in \(M\). This theorem was later quantified in \cite{mns} by Marques, Neves and Song, where they proved the existence of an equidistributed sequence of closed minimal hypersurfaces in \((M,g)\) (for a generic \mt \(g\)). The density and the equidistribution theorems were proved in the Allen-Cahn setting by Gaspar and Guaraco \cite{gg}. Using the multiplicity one theorem, proved by Zhou \cite{zhou2}, Marques and Neves \cite{mn_morse} proved that for a generic (bumpy) \mt \(g\), there exists a sequence of closed, two sided, minimal hypersurfaces \(\{\Si_k\}_{k=1}^{\infty}\) in \((M,g)\) \st Ind\((\Si_k)=k\) and \(\mathcal{H}^n(\Si_k) \sim k^{\frac{1}{n+1}}\). This theorem was previously proved for \(n+1=3\), in the Allen-Cahn setting, by Chodosh and Mantoulidis \cite{cm}.

In \cite{zz1}, Zhou and Zhu developed a min-max theory to construct closed CMC hypersurface in an arbitrary closed Riemannian manifold. More precisely, they proved that if \((M^{n+1},g)\) is a closed Riemannian manifold, \(3\leq n+1 \leq 7\), given any \(c \in \bbr^+\), there exists a closed, \textit{almost embedded} \hy in \((M,g)\) with constant mean curvature \(c\). Here almost embedded means that all the self-intersections are one sided, tangential intersections. (For a precise definition, see Definition \ref{def 3.1}.) Later they generalized this theorem in \cite{zz2} and proved that for a generic function \(h\in C^{\infty}(M)\), there exists a closed, almost embedded \hy in \((M,g)\) with prescribed mean curvature (PMC) \(h\). The min-max construction of the PMC hypersurfaces was used by Zhou \cite{zhou2} to prove the multiplicity one conjecture of Marques and Neves \cite{MN_index}.

Given the above mentioned results regarding the abundance of closed minimal hypersurfaces in closed Riemannian manifolds, it is natural to ask, for a given \(c>0\) and a closed Riemannian manifold \((M,g)\), whether there exist multiple closed \(c\)-CMC hypersurfaces in \((M,g)\). We will show that for any closed Riemannian manifold \((M,g)\), the number of closed \(c\)-CMC hypersurfaces in \((M,g)\) tends to infinity as \(c \ra 0^+\). This will be a consequence of the following theorem. (We refer to Section \ref{sec.2} and \ref{sec 3} for the relevant notations and definitions.)

\begin{thm}
	\label{main thm}
	Let $(M^{n+1},g)$ be a closed Riemannian manifold, $n+1\geq 3$. Suppose $k \in \bbn$ \st $\omega_{k} < \omega_{k+1}$, $c \in \bbr^{+}$ \st $c\textup{Vol}(g)<\omega_{k+1}-\omega_k$ and \(\et \in \bbr^+\) is arbitrary. Then \te \(\Om \in \cm\) \st \(\del \Om \) is a closed $c$-CMC (\wrt the inward unit normal) \hy with optimal regularity and $\omega_k<\ac(\Omega)<\omega_k + W_0 +\et$.
\end{thm}

As a corollary of Theorem \ref{main thm}, we obtain the following theorem.

\begin{thm}\label{main cor}
	Let $(M^{n+1},g)$ be a closed Riemannian manifold, $n+1\geq 3$. For all $m \in \bbn$, \te a constant $c^{*}(m)>0$, which depends on the \mt \(g\), \st if \(0<c<c^{*}(m)\), \((M,g)\) contains at least \(m\) many closed \(c\)-CMC hypersurfaces with optimal regularity. More quantitatively, there exists a constant \(\ga_0\), depending on \(g\), \st for all \(c>0\), there exist at least \(\ga_0c^{-\frac{1}{n+1}}\) many closed \(c\)-CMC hypersurfaces (with optimal regularity) in \((M,g)\).
\end{thm}

At the end of the proof of Theorem \ref{main cor}, we will obtain an explicit lower bound for \(\ga_0\). 

On the other hand, by the work of Pacard and Xu \cite{px}, if \(\La_M\) is the Lusternik-Shnirelman category of \(M\) (i.e. the minimum number of critical points of a smooth function $f:M\ra \bbr$), for \(c\) sufficiently large, there exist at least \(\La_M\) many closed, embedded \(c\)-CMC hypersurfaces in \((M,g)\).

The proof of Theorem \ref{main thm} relies on the works by Zhou-Zhu \cite{zz1} and Zhou \cite{zhou2}. We will also use the regularity theory of stable CMC hypersurfaces, developed by Bellettini-Wickramasekera \cite{bw1, bw2} and Bellettini-Chodosh-Wickramasekera \cite{bcw}, to extend the theorem of Zhou and Zhu \cite{zz1} in higher dimensions. Finally, Theorem \ref{main cor} is deduced from Theorem \ref{main thm} using the growth estimate of the volume spectrum, proved by Gromov \cite{gro}, Guth \cite{guth}, Marques-Neves \cite{mn_ricci_positive} and Liokumovich-Marques-Neves \cite{lmn}.

\textbf{Acknowledgements.} I am very grateful to my advisor Prof. Fernando Cod\'{a} Marques for many helpful discussions and for his constant support and guidance. I thank him for mentioning the inequality \eqref{e.subadditive property of omega_k} to me. The author is partially supported by NSF grant DMS-1811840.

\section{Notation and Preliminaries}
\label{sec.2}
\subsection{Notation}
Here we summarize the notation which will be frequently used later.

\begin{tabular}{ll}
	$\cm$ & the space of Caccioppoli sets in $M$.\\
	\(\znf\), \(\zn\)  & the space of mod \(2\) flat hypercycles in \(M\), which are boun-\\ &daries, equipped with the \(\cF\) norm and the \(\bF\)-metric.  \\
	\(\om_p=\om_p(M,g)\) & the \(p\)-width of \((M,g)\). \\
	\(W_0=W_0(M,g)\) & the one parameter Almgren-Pitts width of \((M,g)\). \\
	$\cH^s$ & the Hausdorff measure of dimension \(s\). \\
	\(B(p,r)\) & the open ball of radius \(r\), centered at \(p\). \\
	\(A(p,r,R)\) & the open annulus, centered at \(p\), with radii \(r<R\).\\
	\(|\Si|\) & the \(n\)-\vf associated to an \(n\)-rectifiable set \(\Si\).\\
	\(\norm{V}\) & the Radon measure associated to a varifold \(V\).\\
	\(CA\) & the cone over a topological space \(A\). \\
	\(SA\) & the suspension of a topological space $A$. \\
	\(\De^k,[v_0, \dots ,v_k]\) & the standard \(k\)-simplex with vertices \(v_0,\dots, v_k\). \\
	\(C_i(A, \bbz_2)\), \(C^i(A,\bbz_2)\) & the abelian group of \(i\) dimensional singular chains \\ & and cochains in \(A\) with \(\bbz_2\) coefficient. \\
	\(H_i(A,\bbz_2)\), \(H^i(A,\bbz_2)\) & the \(i\)-th singular homology and cohomology group  of \(A\) with \\ & \(\bbz_2\) coefficient.
\end{tabular}	

\subsection{Caccioppoli sets and the $\ac$ functional}
\label{sec.2.2} 
In this subsection, we will briefly recall the notion of the \textit{Caccioppoli set}; further details can be found in Simon's book \cite{sim}. \(E \subset (M,g)\) is called a \textit{Caccioppoli set} if \(E\) is \(\mathcal{H}^{n+1}\)-measurable and \(D\chi_E\), the distributional derivative of the characteristic function \(\chi_E\), is a Radon measure. This is equivalent to
\[\sup \left\{\int_{E}\operatorname{div } \om\; d\mathcal{H}^{n+1}: \om \in \mathfrak{X}^1(M) \text{ and } \|\om\|_{\infty}\leq 1 \right\}<\infty,\] 
where \(\mathfrak{X}^1(M)\) is the space of \(C^1\) vector-fields on \(M\). The space of Caccioppoli sets in $M$ will be denoted by \(\cm\). If \(E \in \cm\), \te an \(n\)-rectifiable set, denoted by \(\del E\), \st the total variation measure \(|D\chi_E|=\mathcal{H}^n\mres\del E\). Further, there exists a \(|D\chi_E|\)-measurable vector-field \(\nu_{\del E}\) \st \(\|\nu_{\del E}\|=1\) \(|D\chi_E|\)-a.e. and 
\[\int_E\text{div } \om\; d\mathcal{H}^{n+1}=\int_{\del E}g\left(\om ,\nu_{\del E}\right),\]
for all $\om \in \mathfrak{X}^1(M)$. \(\nu_{\del E}\) is called the \textit{outward unit normal} to \(\del E\). 
For \(c>0\), the \(\ac\) functional on \(\cm\) is defined by
\[\ac(\Om)=\mathcal{H}^n(\del \Om)-c\mathcal{H}^{n+1}(\Om).\]
Let \(X\) be a vector-field on \(M\) and \(\vp_t\) be the flow of \(X\). The first variation of \(\ac\) along \(X\) is given by
\[\de \ac\big|_{\Om}(X)=\frac{d}{dt}\bigg|_{t=0}\ac(\vp_t(\Om))=\int_{\del \Om}\operatorname{div}_{\del \Om}X d\mathcal{H}^n-c\int_{\del \Om}g\left(X,\nu_{\del \Om}\right)d\mathcal{H}^n.\]
Suppose \(\Om\) is an open set \st \(\del \Om\) is a smooth hypersurface with mean curvature vector \(H\). Then \(\text{div}_{\del \Om}X = -g(X,H)\). \tf \(\Om\) is a critical point of the \(\ac\) functional if and only if $H=-c\nu_{\del\Om}$, i.e. \(\del \Om\) has constant mean curvature \(c\) \wrt the inward unit normal \(-\nu_{\del\Om}\). In this case, the second variation of \(\ac\) along \(X\) (which is assumed to be normal along \(\del \Om\)) is given by
\begin{align*}
\label{e.second variation}
\de^2 \ac\big|_{\Om}(X,X) &=\frac{d^2}{dt^2}\bigg|_{t=0}\ac(\vp_t(\Om))\nonumber\\ &
=\int_{\del \Om}\left(|\nabla^{\perp}X|^2-\text{Ric}(X,X)-\left|A_{\del\Om}\right|^2|X|^2\right)d\mathcal{H}^n.
\end{align*}
Here \(A_{\del\Om}\) stands for the second fundamental form of \(\del \Om\).

\subsection{The space of hypercycles and the volume spectrum}
\label{sec.2.3}
Let us introduce the following notation. \(\mathbf{I}_l(M^{n+1};\bbz_2)\) is the space of \(l\)-dimensional mod \(2\) flat chains in \(M\); we only need to consider \(l=n,n+1\). \(\mathcal{Z}_n(M;\bbz_2)\) denotes the space of flat chains \(T \in \mathbf{I}_n(M;\bbz_2)\) \st \(T=\del U\) for some \(U \in \mathbf{I}_{n+1}(M;\bbz_2)\). For \(T \in \mathbf{I}_n(M;\bbz_2)\), \(\md{T}\) stands for the varifold associated to \(T\) and \(\norm{T}\) is the Radon measure associated to \(\md{T}\). \(\cF\) and \(\bM\) denote the flat norm and the mass norm on \(\mathbf{I}_l(M;\bbz_2)\). When \(l=n+1\), these two norms coincide. The \(\bF\)-metric on the space of currents is defined as follows.
\begin{align*}
&\bF(U_1,U_2)=\cF(U_1-U_2)+\bF(\md{\del U_1},\md{\del U_2})\; \text{ if } U_1, U_2 \in \mathbf{I}_{n+1}(M;\bbz_2);\\
&\bF(T_1,T_2)=\cF(T_1-T_2)+\bF(\md{T_1},\md{ T_2})\; \text{ if } T_1, T_2 \in \mathcal{Z}_{n}(M;\bbz_2).
\end{align*}
\(\znf\) and \(\zn\) will stand for the space \(\cZ_n(M;\bbz_2)\) equipped with the \(\cF\) norm and the \(\bF\)-metric respectively.
We will identify \(\cm\) with \(\mathbf{I}_{n+1}(M;\bbz_2)\), i.e. \(E \in \cm\) will be identified with \(\db{E}\), the current associated with \(E\). Similarly, \(\del E\) (with \(E \in \cm\)) will be identified with \(\db{\del E}=\del \db{E}\).

By the constancy theorem, if \(\Om_1, \Om_2 \in \cm\) \st \(\del \Om_1 = \del \Om_2\), then either \(\Om_1=\Om_2\) or \(\Om_1=M-\Om_2\). In \cite{mn_morse}, Marques and Neves proved that the space \((\cm, \cF)\) is contractible and the boundary map
\[\del:(\cm, \cF) \ra \mathcal{Z}_n(M;\cF;\bbz_2)\]
is a \(2\)-sheeted covering map. By the definition of the \(\bF\)-metric on \(\cm\) and \(\mathcal{Z}_n(M;\bbz_2)\), this implies
\[\del:(\cm, \bF) \ra \zn\]
is also a \(2\)-sheeted covering map. Furthermore, \(\pi_1(\znf)=\bbz_2\) and for \(i\geq 2\), \(\pi_i(\znf)=0\). It was also proved in \cite{mn_morse} that \(\znf\) is weakly homotopy equivalent to \(\mathbb{RP}^{\infty}\). Let \(\overline{\la}\) denote the unique nonzero element of \(H^1\left(\znf,\bbz_2\right)(=\bbz_2)\); then the cohomology ring \(H^*\left(\znf,\bbz_2\right)=\bbz_2[\overline{\la}]\).

\(X\) is called a \textit{cubical complex} if \(X\) is a subcomplex of \([0,1]^l\) for some \(l \in \bbn\). The cells of \([0,1]^l\) are of the form $\al_1\otimes\al_2\otimes\dots \otimes\al_l$, where each $\al_i\in \{[0],[1],[0,1]\}$. By \cite{bp}*{Chapter 4}, every simplicial complex is homeomorphic to a cubical complex and vice-versa. 

Let $X$ be a cubical complex. A \cts map \(\Ph:X\ra \znf\) is called a \textit{\(k\)-sweepout} if 
\[(\vp^*\overline{\la})^k\neq 0 \in H^k(X,\bbz_2).\] 
\(\Ph:X\ra \mathcal{Z}_n(M;\bbz_2)\) is said to have the property {\it no concentration of mass} if
\[\lim_{r \ra 0}\sup \{\norm{\Ph(x)}(B(p,r)):x \in X, p \in M\}=0.\]
A continuous map \(\Ph:X \ra \zn\) has no concentration of mass \cite{MN_index}*{Proof of Theorem 3.8}. The set of all the \(k\)-sweepouts with no concentration of mass is denoted by \(\mathcal{P}_k\). The \textit{\(k\)-width} of \((M,g)\) is defined by
\[\om_{k}(M,g)=\inf_{\Ph \in \mathcal{P}_k}\sup\left\{\bM(\Ph(x)):x \in \text{dmn}(\Ph)\right\},\]
where \(\text{dmn}(\Ph)\) is the domain of \(\Ph\). The volume spectrum \(\{\om_k\}_{k=1}^{\infty}\) satisfies the following inequality, proved by Gromov \cite{gro}, Guth \cite{guth}, Marques-Neves \cite{mn_ricci_positive} and Liokumovich-Marques-Neves \cite{lmn}.
\begin{equation}
\label{e. sublinear growth}
\om_{k} \geq K_0k^{\frac{1}{n+1}}\quad \forall\; k \in \bbn,
\end{equation}
where \(K_0>0\) is a constant, which depends on the \mt \(g\).

Let \(\mathcal{S}\) be the set of all continuous maps \(\La:[0,1] \ra (\cm, \cF)\) \st \(\La(0)=M,\;\La(1)=\emptyset\) and \(\del \circ \La :[0,1] \ra \znf\) has no concentration of mass. The \textit{one parameter Almgren-Pitts width of \((M,g)\)}, denoted by \(W_0(M,g)\), is defined by 
\[W_0(M,g)=\inf_{\La \in \mathcal{S}}\sup\left\{\bM(\del \La(t)):t \in [0,1]\right\}.\]

\subsection{Some notions from topology}
\label{sec.2.4}
Given a topological space \(A\), let \(CA\) denote the \textit{cone over \(A\)}, which is defined as follows.
\[CA=\frac{A\times[0,1]}{\sim},\]
where the equivalence relation `\(\sim\)' collapses \(A \times \{1\}\) to a point. The cone construction is functorial, i.e. if \(f:A\ra B\) is a \cts map between the topological spaces \(A\) and \(B\), there exists a \cts map \(Cf:CA \ra CB\) defined by \(Cf\left(a,t\right)=\left(f(a),t\right)\). Here, by the abuse of notation, we are denoting an element of \(CA\) by a pair \((a,t)\) with \(a \in A\) and \(t \in [0,1]\). We note that if \(p_0\in CA\) is the collapsed image of \(A \times \{1\}\) and \(q_0 \in CB\) is the collapsed image of \(B \times \{1\}\), then \(Cf(p_0)=q_0\).

The \textit{suspension of \(A\)}, denoted by \(SA\), is defined as follows.
\[SA=\frac{A\times[-1,1]}{\sim},\]
where the equivalence relation `\(\sim\)' collapses \(A \times \{1\}\) to a point and \(A \times \{-1\}\) to another point. \(SA=C_+A\cup C_-A\), where
\[C_+A= \frac{A\times[0,1]}{\sim},\quad \quad C_-A= \frac{A\times[-1,0]}{\sim}.\]
\(C_+A\) and \(C_-A\) are cones over \(A\). Since the cone over a topological space is contractible, using van Kampen's Theorem, one can show that if \(A\) is path-connected, \(SA\) is simply-connected. We note that the cone over a \(k\)-simplex is a \((k+1)\)-simplex. \tf if \(A\) is a cubical complex, \(CA\) and \(SA\) are also cubical complexes.

\section{Regularity of the min-max CMC hypersurfaces in higher dimensions}\label{sec 3}
In this section, we will modify the argument of Zhou and Zhu \cites{zz1, zhou2} to extend their min-max theory (Theorem \ref{min-max thm}) in higher dimensions. We will use the compactness theory of stable CMC \hys (Theorem \ref{cptness thm}), developed by Bellettini-Wickramasekera \cite{bw1, bw2} and Bellettini-Chodosh-Wickramasekera \cite{bcw}. When the ambient dimension \(3\leq n+1\leq 7\), the compactness theorem for stable CMC \hys was also proved by Zhou and Zhu \cite{zz1}.

\begin{defn}\label{def 3.1}
	\(\Si\subset (U^{n+1},g)\) is called a \textit{\(c\)-CMC \hy \w} (\(c\geq 0\)) if $\Si$ is a closed \(n\)-rectifiable subset of \((U,g)\) and there exists a closed subset of \(\Si\), denoted by \(sing(\Si)\), \st the following conditions are satisfied.
	\begin{itemize}
		\item[(i)]  \(\cH^s(sing(\Si))=0\) for all \(s>n-7\).
		\item[(ii)] \(reg(\Si)=\Si \setminus sing(\Si)\) is {\it almost embedded}, i.e. if \(p\in reg(\Si)\) \st \(\Si\) is not embedded near \(p\), \te an open set \(B \subset U\), containing \(p\), \st
		
		\textbullet\- \(B \cap \Si=\Si_1 \cup \Si_2\), where \(\Si_1\) and \(\Si_2\) are smoothly embedded hypersurfaces in \(B\) (since \(\Si\) is not embedded at \(p\), we necessarily have \(p \in \Si_1 \cap \Si_2\));
		
		\textbullet\- \(B \setminus \Si_1=B_1 \cup B_2\), where \(B_1\) and \(B_2\) are connected and \(\Si_2 \subset \overline{B}_1\) (hence \(\Si_1\) and \(\Si_2\) intersect each other tangentially).
			
		\item[(iii)] $|\Si|$ has bounded first variation and \(reg(\Si)\) has constant mean curvature \(c\) \wrt a globally defined unit normal \(\nu\). (If \(c=0\), we do not require that \(reg(\Si)\) is two sided.)
	\end{itemize}
\end{defn}
Given such a \(\Si\), following \cite{zz1}, let us define
\[\cR(\Si)=\left\{p \in reg(\Si): \Si \text{ is embedded near } p \right\};\;  \; \cS(\Si) = reg(\Si)\setminus \cR(\Si).\]
If \(\Si\) is a minimal \hy (i.e. $c=0$) with optimal regularity, then, by the maximum principle, \(\cS(\Si)=\emptyset\). If $c>0$, by \cite{zz1}*{Proposition 2.9} and \cite{bw1}*{Remark 2.6}, \(\cS(\Si)\) is \((n-1)\)-rectifiable.

\begin{defn}
	Let \(\Si \subset U\) be a \(c\)-CMC hypersurface with optimal regularity (\(c\geq 0\)). \(\Si\) is called \textit{stable} if for all compactly supported, normal vector-field \(\scrX\) on \(reg(\Si)\),
	\[\int_{reg(\Si)}\left(|\nabla^{\perp}\scrX|^2-\text{Ric}(\scrX,\scrX)-\left|A_{\Si}\right|^2|\scrX|^2\right)d\mathcal{H}^n\geq 0.\]
Here \(A_{\Si}\) is the second fundamental form of \(\Si\), which is defined on \(reg(\Si)\).
\end{defn}

\begin{thm}[\cite{bw1,bw2,bcw}, \cite{zz1} when \(3\leq n+1\leq 7\)]\label{cptness thm}
	Suppose \(\Si_k \subset U \) is a \seq of stable, \(c_k\)-CMC hypersurfaces with optimal regularity, \(\sup_k \cH^n(\Si_k)<\infty\) and \(c_k \ra c_{\infty}\). Then, possibly after passing to a subsequence, \(|\Si_k|\) converges to an integral \(n\)-\vf \(W\) \st \(\spt\|W\|=\Si_{\infty}\) is a closed \(c_{\infty}\)-CMC \hy with optimal regularity. The convergence is smooth on compact subsets of \(reg(\Si_{\infty})\). Moreover, if each \(\Si_k=\del \Om_k\) for some \(\Om_k \in \mathcal{C}(U)\) and \(c_{\infty}>0\), then \(\Si_{\infty}=\del \Om_{\infty}\) for some \(\Om_{\infty}\in \mathcal{C}(U)\) and the density of \(W\) is \(1\) on \(\cR(\Si_{\infty})\) and \(2\) on \(\cS(\Si_{\infty})\).
\end{thm}

The following definition is taken from \cite{zhou2}. We remark that in \cite{zhou2}*{Definition 1.1}, the homotopies are continuous in the flat norm. However, if we examine the proof of Theorem 1.7 in \cite{zhou2}, we see that the homotopies can be taken to be \cts in the \(\bF\)-metric. (In \cite{zhou2}*{Proposition 1.14 and 1.15} it is proved that the Almgren extensions are homotopic to each other and to the initial map in the \(\bF\)-metric; see also \cite{mn_morse}*{Section 3}.)

\begin{defn}
	 \sps $X$ is a cubical complex, $Z\subset X$ is a subcomplex, $F_0 : X \rightarrow (\cm,\bF)$ is a continuous map. Let $\Pi$ denote the set of all sequence of \cts maps \(\left\lbrace F_i:X \rightarrow (\cm,\bF)\right\rbrace_{i=1}^{\infty}\) \st for every \(i \in \bbn\), \te homotopy \(G_i:X \times [0,1] \rightarrow (\cm,\bF)\) with the properties \(G_i(-,0)=F_0\), \(G_i(-,1)=F_i\) and 
	 \begin{equation*}
	 	\limsup_{i \rightarrow \infty}\sup \left\lbrace \bF(G_i(x,s), F_0(x)):x \in Z, s\in [0,1]\right \rbrace=0.
	 \end{equation*}
	 \(\Pi\) is called the \textit{\((X,Z)\)-homotopy class of \(F_0\)}. For \(c>0\), we define
	 \[\bL^c(\Pi)=\inf_{\{F_i\} \in \Pi}\limsup_{i \rightarrow \infty}\sup_{x \in X}\left\lbrace\ac(F_i(x))\right\rbrace.\]
\end{defn}

To prove the following theorem, we only need to slightly modify the argument given in the paper \cite{zz1}. We will only mention the necessary modifications; further details can be found in \cite{zz1}*{Section 5 and 6}.

\begin{thm}[\cite{zz1,zhou2}]\label{min-max thm}
	Let \((M^{n+1},g)\) be a closed Riemannian manifold, \(n+1 \geq 3\), and $c>0$. \sps  $F_0 : X \rightarrow (\cm,\bF)$ is a continuous map, $Z \subset X$ and $\Pi$ is the \((X,Z)\)-homotopy class of \(F_0\) \st \[\bL^c(\Pi)> \sup_{x \in Z}\ac(F_0(x)).\] Then \te \(\Om \in \cm\) \st $\ac(\Om)=\bL^c(\Pi)$ and $\del \Om$ is a closed $c$-CMC (\wrt the inward unit normal) \hy with optimal regularity.
\end{thm}

\begin{proof}
	Throughout the proof, we will use the notation used in \cite{zz1}. It suffices to show that if \(V \in \mathcal{V}_n(M)\) has \(c\)-bounded first variation in \(M\) and is \(c\)-almost minimizing in small annuli, then \(V=|\Si|\), where \(\Si\) is a \(c\)-CMC \hy with optimal regularity. We will assume that \(M\) is isometrically embedded inside some Euclidean space \(\bbr^L\). In the rest of this section, \(B(p,r)\) and \(A(p,r,R)\) will respectively denote the ball and the annulus in \(\bbr^L\) with respect to the Euclidean metric.
	
    For \(p \in \bbr^L\) and \(r>0\), \(\bet_{p,r}:\bbr^L\ra \bbr^L\) is defined by \(\bet_{p,r}(x)=\frac{x-p}{r}\). We start with the following proposition.

\begin{pro}\label{blow up is regular}
	Let \(V\) be \(c\)-almost minimizing in \(U\). Suppose \(\{p_i\}\) is a sequence in \(U\) converging to \(p \in U\) and \(r_i>0\) is a sequence converging to \(0\). If \(\overline{V}\) is the varifold limit of \((\bet_{p_i,r_i})_{\#}V\), then \(\overline{V}\) is induced by a minimal \hy with optimal regularity (possibly with integer multiplicities). In particular, if \(V \in \mathcal{V}_n(M)\) has \(c\)-bounded first variation in \(M\) and is \(c\)-almost minimizing in small annuli, then \(V\) is integer rectifiable and every \(C \in \vt(V,p)\) is induced by a minimal \hy with optimal regularity  (possibly with integer multiplicities).
\end{pro}

The proof of this proposition is same as the proof of Lemma 5.10 and Proposition 5.11 in \cite{zz1}. One can show that \(\overline{V}\) satisfies the good replacement property. This implies the proposition by the results of Schoen-Simon \cite{ss} and De Lellis-Tasnady \cite{dt}.

\begin{lem}\label{half space theorem}
	Let \(H=\{x\in\bbr^{n+1}:x_{n+1}=0\}\) and \(S\neq H\) be another hyperplane in \(\bbr^{n+1}\), which is not parallel to \(H\). Suppose \(W\) is a stationary, integral \vf in \(\bbr^{n+1}\) \st \(\spt\|W\|\) is a minimal \hy \w and \(W\mres\{x\in \bbr^{n+1}:x_{n+1}>0\}=m|S|\mres\{x\in \bbr^{n+1}:x_{n+1}>0\}\). Then \(W=m|S|\).	
\end{lem}
\begin{proof}
	The proof is an adaptation of the argument of Hoffman and Meeks \cite{hoff-meeks}*{Proof of Theorem 3} in the singular setting. Let \(\Xi=\spt \norm{W}\). By our assumption, 
	\begin{equation}
	\cH^s(sing(\Xi))=0 \; \text{ for all } s>n-7. \label{e. smallness of sing Sigma}
	\end{equation}
	\tf \te a collection \(\{\Xi_j\}_{j=0}^{J}\) of disjoint, connected minimal \hys \w and positive integers $\{m_j\}_{j=0}^J$ \st
	\[W=\sum_{j=0}^J m_j|\Xi_j|.\]
	Since \(\Xi \cap \{x_{n+1}>0\}=S\cap\{x_{n+1}>0\}\), by \eqref{e. smallness of sing Sigma}, we have \(S \subset reg(\Xi)\) (see \cite{ss}*{(7.23)}). By the theorem of Ilmanen \cite{Ilmanen}*{Theorem A(ii)} and Wickramasekera \cite{W}, \eqref{e. smallness of sing Sigma} implies that the connected components of \(reg(\Xi)\) are precisely \(\{reg(\Xi_j)\}_{j=0}^{J}\). Thus \(S=reg(\Xi_l)\), for some \(l \in \{0,\dots,J\}\). Let us assume that \(S= reg(\Xi_0)\); \eqref{e. smallness of sing Sigma} gives \(S=\Xi_0\) and hence \(m_0=m\). If \(W\neq m|S|\), let us consider \(\Xi_1\). Since \(W\mres\{x_{n+1}>0\}=m|S|\mres\{x_{n+1}>0\}\), \(\Xi_1 \subset \{x_{n+1}\leq 0\}\). We claim that \(\Xi_1 \subset \{x_{n+1}< 0\}\). Indeed, by the maximum principle proved by White \cite{white}*{Theorem 4} and Solomon-White \cite{sw}, if \(\Xi_1 \cap H \neq \emptyset\), then \(H \subset \Xi_1\), which contradicts the fact that \(\Xi_1\) is disjoint from \(\Xi_0=S\). As $\Xi_1$ is connected, \(\Xi_1\) lies in a region \(\mathscr{R}\) which is the intersection of the open half-space \(\{x_{n+1}<0\}\) and another open half-space defined by \(S\). Using an orthogonal transformation, we can obtain a new coordinate system \(\{y_1,\dots,y_{n+1}\}\) \st \(\mathscr{R}\subset \{y_{n+1}>0\}\), \(\del \mathscr{R}\) is a graph over the hyperplane $\{y_{n+1}=0\}$ and the \(y_1\)-axis is contained in \(H\cap S\). As \(\Xi_1\) is a closed set and is disjoint from \(H \cup S\), \te \(r>0\) \st 
	\[B_r=\left\{y \in \bbr^{n+1}:(y_1-r)^2+y_2^2+\dots+y_{n+1}^2\leq r^2\right\}\]
	is disjoint from \(\Xi_1\). We note that \(\del B_r \cap \del \mathscr{R}\) is a graph over the boundary of a convex domain in \(\{y_{n+1}=0\}\). \tf by \cite{gt}*{Theorem 16.8}, \te a smooth minimal \hy \(N\) with boundary \(\del N = \del B_r \cap \del \mathscr{R}\). By the convex hull property, \(N \subset B_r\); hence \(N \cap \Xi_1=\emptyset\). For \(t \in [1,\infty)\), let us consider
	\[N_t=\{tp:p \in N\}.\]
	We claim that for all \(t\in [1,\infty)\), \(N_t \cap \Xi_1 =\emptyset\). Otherwise, \te a smallest \(\tau>1\) \st \(N_{\tau} \cap \Xi_1 \neq\emptyset\). This implies, by \cite{sw}*{Step 1, p. 687}, \(N_{\tau} \subset \Xi_1\), which contradicts the fact that \(\Xi_1 \cap \del \mathscr{R}=\emptyset\). However, 
	\[\mathscr{R}\cap\{y_1>0\} \subset B_r \cup \left(\cup_{t\in [1,\infty)}N_t\right).\]
	\tf \(\Xi_1 \subset\mathscr{R}\cap\{y_1\leq 0\} \).
	An analogous argument gives that for any \(T \in \bbr\), \(\Xi_1 \subset\mathscr{R}\cap\{y_1\leq T\} \), which forces \(\Xi_1=\emptyset\). This finishes the proof of the lemma.
\end{proof}

We proceed as in \cite{zz1}. Let us fix \(p\in \spt\norm{V}\). We choose \(0<r_0<r_{am}(p)\) \st 
\begin{itemize}
	\item for all \(0<r\leq r_0\), \(\del B(p,r)\cap M\) is a smooth \hy in \(M\) with mean curvature greater than \(c\);
	\item \(x\mapsto d_{\bbr^L}(p,x)\) is a smooth function on \((B(p,r_0)\cap M)\setminus\{p\}\).
\end{itemize}
For \(r \leq r_0\) and \(W \in \mathcal{V}_n(M)\), if \(W \neq 0\) in \(B(p,r)\) and has \(c\)-bounded first variation in \(B(p,r)\cap M\),  we have
\begin{equation}
\emptyset\neq \spt \norm{W}\cap \del B(p,r)=\overline{\spt \norm{W}\setminus \overline{B}(p,r)}\cap \del B(p,r);
\label{e. approximation from outside}
\end{equation}
\begin{equation}
\spt \norm{W}\cap \overline{B}(p,r)=\overline{\bigcup_{0<s<r}\spt (\norm{W}\mres B(p,s))\cap \del B(p,s)}.
\label{e. approximation from inside}
\end{equation}
Equation \eqref{e. approximation from inside} is a consequence of \cite{zz1}*{Lemma 6.2} (see \cite{dt}*{equation (5.6) and A.3. Proof of Lemma 5.4}). Let us fix \(0<s<t<r_0\). \(V^*\) is a \(c\)-replacement of \(V\) in \(\cl (A(p,s,t)\cap M)\). By the regularity of the minimizers for the \(\ac\) functional \cite{zz1}*{Theorem 2.14; \cite{M}} and Theorem \ref{cptness thm}, there exists \(\Si_1\), a \(c\)-CMC \hy in \(A(p,s,t)\cap M\) with optimal regularity, \st \(\Si_1=\del \Om_1 \cap A(p,s,t)\cap M\) for some \(\Om_1 \in \cm\) and
\[V^{*}\mres(A(p,s,t)\cap M)=|\Si_1|.\]
\(\cS(\Si_1)\) is contained in \(\cup_k \Si^{(k)}_1\), a countable union of \((n-1)\)-dimensional submanifolds. We choose \(s_2\in (s,t)\) \st \(\del B(p,s_2)\) intersects \(reg(\Si_1)\) and all the \(\Si^{(k)}_1\)'s transversally. For \(s_1 \in (0,s)\), let us consider \(c\)-replacement \(V^{**}\) of \(V^*\) in \(\cl (A(p,s_1,s_2)\cap M)\). There exists \(\Si_2\), a \(c\)-CMC \hy in \(A(p,s_1,s_2)\cap M\) with optimal regularity, \st \(\Si_2=\del \Om_2 \cap A(p,s_1,s_2)\cap M\) for some \(\Om_2 \in \cm\) and
\[V^{**}\mres(A(p,s_1,s_2)\cap M)=|\Si_2|.\]
We need to show that \(\Si_2\) smoothly glues with \(\Si_1\) across \(\del (B(p,s_2) \cap M)\).
\begin{lem}
	\label{tangent space matching} For all \(x \in reg(\Si_1)\cap \del (B(p,s_2)\cap M)\), \(\vt(V^{**},x)=\{\Theta^n(\norm{V^*},x)|T_x\Si_1|\}\), where \(T_x\Si_1\) denotes the tangent space of \(\Si_1\) at \(x\) with multiplicity one.
\end{lem}
\begin{proof}
	Let \(C\in\vt(V^{**},x)\). Denoting \(\La=\del (B(p,s_2) \cap M)\), since \(V^{**}=V^*=|\Si_1|\) in \(A(p,s_2,t)\cap M\), we must have \(C = \Theta^n(\norm{V^*},x)|T_x\Si_1|\) in an open half-space of \(T_xM\), defined by \(T_x\La\). (By the transversality assumption, \(T_x\Si_1 \neq T_x\La\).) Using Proposition \ref{blow up is regular} and Lemma \ref{half space theorem}, we obtain \(C=\Theta^n(\norm{V^*},x)|T_x\Si_1|\).
\end{proof}
As argued in \cite{zz1}*{Step 2}, using \eqref{e. approximation from outside}, one can show that
\begin{equation*}
\label{eq 3.1}\cl(\Si_2)\cap \del (B(p,s_2) \cap M) \subset \Si_1\cap \del (B(p,s_2) \cap M).
\end{equation*}
Moreover, smallness of \(sing(\Si_1)\) and Lemma \ref{tangent space matching} imply that
\[\Si_1\cap \del (B(p,s_2) \cap M)=\cl (reg(\Si_1)\cap \del (B(p,s_2) \cap M))\subset \cl(\Si_2)\cap \del (B(p,s_2) \cap M).\]
Hence,
\[\cl(\Si_2) \cap \del (B(p,s_2)\cap M)= \Si_1 \cap \del (B(p,s_2)\cap M),\] 
i.e. \(\Si_1,\; \Si_2\) glue continuously along \(\del B(p,s_2)\). Furthermore,  Lemma \ref{tangent space matching}, \cite{sim}*{Theorem 3.2(2)} and smallness of \(sing(\Si_1)\) also give that
\(\norm{V^{**}}(\del B(p,s_2))=0,\)
which implies, by \cite{zz1}*{Claim 1(c) in Section 6},
\[V^{**}\mres (A(p,s_1,t)\cap M)=\md{\del \Om^{**}}\mres (A(p,s_1,t)\cap M),\]
for some \(\Om^{**}\in \cm.\)

Arguing as in \cite{zz1}*{Step 2}, we can also show that \(reg(\Si_1)\) and \(reg(\Si_2)\) glue smoothly across \(\del( B(p,s_2)\cap M)\); we need to use Lemma \ref{half space theorem} in place of the half space theorem of \cite{hoff-meeks} and Theorem \ref{cptness thm} in place of \cite{zz1}*{Theorem 2.11}. Let us give some details. Let \(\nu_1\), \(\nu_2\) denote the outward unit normal of \(reg(\Si_1)\) and \(reg(\Si_2)\) respectively;
\[\Ga=reg(\Si_1)\cap \del (B(p,s_2)\cap M),\quad \cR(\Ga)=\cR(\Si_1)\cap \Ga, \quad \cS(\Ga)=\cS(\Si_1)\cap \Ga.\]
Suppose \(x \in \cR(\Ga)\), \(y_i \in reg(\Si_2)\) converges to \(x\). It is possible to choose \(x_i \in \cR(\Ga)\) \st \(r_i=\md{x_i-y_i}\) converges to \(0\) and hence \(x_i \ra x.\) Using the fact that \(V^{**}\mres (A(p,s_2,t)\cap M)=\md{\Si_1}\mres (A(p,s_2,t)\cap M)\), Proposition \ref{blow up is regular} and Lemma \ref{half space theorem}, one can prove that 
\[\lim_{i \ra \infty}(\et_{x_i,r_i})_{\#}V^{**}=\md{T_x\Si_1}\quad \text{in the sense of varifolds.}\]
This, along with Theorem \ref{cptness thm}, implies that \(\et_{x_i,r_i}(\Si_2\cap B(y_i, r_i/2))\) smoothly converges to a domain in \(T_x\Si_1\), which gives
\begin{equation}\label{C1.1 convergence}
\lim_{y \ra x,\; y\in reg(\Si_2)}\nu_2(y)=\nu_1(x)\; \text{ uniformly in }x\text{ on compact subsets of }\cR(\Ga).
\end{equation}

Similarly, if \(x \in \cS(\Ga)\) and \(y_i \in reg(\Si_2)\) converges to \(x\), it is possible to choose \(x_i \in \Ga\) \st \(r_i=\md{x_i-y_i}\) converges to \(0\) and hence \(x_i \ra x.\) Using the fact that \(V^{**}\mres (A(p,s_2,t)\cap M)=\md{\Si_1}\mres (A(p,s_2,t)\cap M)\), Proposition \ref{blow up is regular} and Lemma \ref{half space theorem}, one can show that 
\begin{equation*}
\lim_{i \ra \infty}(\et_{x_i,r_i})_{\#}V^{**}=
\begin{cases}
\md{T_x\Si_1}+\md{\ta_v T_x \Si_1} & \text{ if } \liminf_{i \ra \infty} d_{\bbr^L}(x_i, \cS(\Ga))/r_i=\infty;\\
2 \md{T_x\Si_1} & \text{ if } \liminf_{i \ra \infty} d_{\bbr^L}(x_i, \cS(\Ga))/r_i<\infty,
\end{cases}
\end{equation*}
where \(T_x\Si_1\) denotes the tangent space of \(\Si_1\) at \(x\) with multiplicity one and \(\ta_v\) denotes translation by \(v \in (T_x\Si_1)^{\perp}\subset T_xM\) (\(v\) might be \(\infty\), in which case \(\ta_v T_x \Si_1=\emptyset\)). As before, this, together with Theorem \ref{cptness thm}, gives that
\begin{equation}\label{C1.2 convergence}
\lim_{y \ra x,\; y\in reg(\Si_2)}[T_y\Si_2]=[T_x\Si_1]\; \text{ uniformly in }x\text{ on compact subsets of } \cS(\Ga),
\end{equation}
where \([T_y\Si_2]\) and \([T_x\Si_1]\) respectively denote the un-oriented, multiplicity one tangent spaces of \(\Si_2\) and \(\Si_1\). \eqref{C1.1 convergence}, \eqref{C1.2 convergence} and elliptic regularity imply that \(reg(\Si_1)\) and \(reg(\Si_2)\) glue smoothly across \(\del( B(p,s_2)\cap M)\). Hence, by unique continuation and by the smallness of \(sing(\Si_1)\) and \(sing(\Si_2)\), \(\Si_2=\Si_1\) in \(A(p,s,s_2)\cap M\).

Since we will vary \( s_1 \in (0,s)\), let us use the notation \(V^{**}_{s_1} \text{ and }\Si_{s_1}\) instead of \(V^{**} \text{ and }\Si_2\). As argued in \cite{zz1}*{Step 3}, the above discussion implies that
\[\Si=\{p\}\cup \bigcup_{0<s_1<s}\Si_{s_1}\]
is a \(c\)-CMC \hy \w in \(B(p,s_2)\) (if necessary, \(p\) can be absorbed inside \(sing(\Si)\)). Further, for any \(s_1<s\),
\begin{equation}
	V^{**}_{s_1}\mres (A(p,s_1,s_2)\cap M)=|\Si|\mres (A(p,s_1,s_2)\cap M). \label{e.V**_s1 matches with V**_s1'}
\end{equation}

Following the argument of Schoen-Simon \cite{ss} and De Lellis-Tasnady \cite{dt}, we will show that \(\spt \norm{V} \cap B(p,s)=\Si \cap B(p,s)\). By the definition of \(c\)-replacement, for all \(0<s_1'\leq s_1\),
\[V^{**}_{s_1}\mres (B(p,s_1')\cap M)=V^{*}\mres (B(p,s_1')\cap M)=V\mres (B(p,s_1')\cap M).\]
Hence, using \eqref{e. approximation from inside},
\begin{equation}
\spt \norm{V}\cap \del (B(p,s_1)\cap M) = \spt \big\|V^{**}_{s_1}\big\|\cap \del (B(p,s_1)\cap M). \label{e. V matches with V^**s_1'}
\end{equation}
Moreover, by \eqref{e.V**_s1 matches with V**_s1'} and \eqref{e. approximation from outside},
\begin{equation}\label{e. V matches with V^**s_1}
\spt \big\|V^{**}_{s_1}\big\|\cap \del (B(p,s_1)\cap M)=\Si \cap \del (B(p,s_1)\cap M).
\end{equation}
\tf combining \eqref{e. V matches with V^**s_1'} and \eqref{e. V matches with V^**s_1},
\[\spt \norm{V}\cap \del (B(p,s_1)\cap M)= \Si \cap \del (B(p,s_1)\cap M),\]
for all \(0<s_1<s\). By varying \(s_1 \in (0, s)\), we obtain \(\spt \norm{V} \cap B(p,s)=\Si \cap B(p,s)\). (We assumed that \(p \in \spt \norm{V}\) and by our definition of \(\Si\), \(p \in \Si\).)

As \(p \in \spt \norm{V}\) is arbitrary, we have proved that \(\spt \norm{V}\) is a \(c\)-CMC \hy with optimal regularity. By a slight abuse of notation, let us denote \(\spt \norm{V}\) by \(\Si\). We will show that if \(p\in \cR(\Si)\), \te \(r>0\) \st  \(\Th^n(\norm{V},y)=1\) for all \(y \in (\cR(\Si) \cap B(p,r))\setminus \{p\}\). This will imply that \(V =|\Si|\) as \(\cH^n(\Si \setminus \cR(\Si))=0\).

The argument is similar to \cite{zz1}*{Proof of Claim 6}. Since \(p \in \cR(\Si)\), one can choose \(0<r<r_0\) (where \(r_0\) is as chosen after the proof of Lemma \ref{half space theorem}), \st \(B(p,r)\cap \Si \subset \cR(\Si)\) and for all \(0<r'\leq r\), \(\del B(p,r')\) intersects \(\Si\) transversally. Let \(\rh = |p-y|\) and we consider the second replacement \(V^{**}_{\rh}\). Since \(V^{**}_{\rh}\mres (A(p,\rh,s)\cap M)=|\Si|\mres (A(p,\rh,s)\cap M)\) for some \(s>\rh\), by the transversality assumption, Proposition \ref{blow up is regular} and Lemma \ref{half space theorem}, \(\vt(V^{**}_{\rh},y)=\{|T_y\Si|\}\). As \(V \mres (B(p,\rh)\cap M)=V^{**}_{\rh}\mres (B(p,\rh)\cap M)\), one can again use Proposition \ref{blow up is regular} and Lemma \ref{half space theorem} to conclude that \(\vt(V,y)=\{|T_y\Si|\}\). Hence \(\Th^n(\norm{V},y)=1\).
\end{proof}

\section{Proof of the main theorems}
\vspace{1\baselineskip}
\begin{proof}[Proof of Theorem \ref{main thm}]
	The theorem will be proved in five parts.
	
	\textbf{Part 1.} We choose \(\de\) \st 
	\begin{equation}\label{E:defn of delta}
	0<\de<\om_{k+1}-\om_k-c\vol(g)\quad \text{ and }\quad \de<\et/2.
	\end{equation} 
	By the definition of $\om_k$, \te a map \(\Ph:X \ra \znf \) with no concentration of mass, \st \(X\) is a connected cubical complex, \(\Ph\) is a $k$-sweepout and 
	\begin{equation}
	\sup_{x \in X}\bM(\Ph(x)) \leq \om_k + \de/2.\label{e. mass of Phi in X}
	\end{equation}
	Following the argument in \cite{zhou2}, $\Ph$ is a \(k\)-sweepout implies that 
	\[\Ph^*:H^1(\znf,\bbz_2)\ra H^1(X,\bbz_2)\]
	is non-zero; hence,
	\[\Ph_{*}:\pi_1(X)\ra \pi_1(\znf)(=\mathbb{Z}_2)\]
	is onto. Thus, \(\ker(\Ph_*)\) is an index \(2\) subgroup of \(\pi_1(X)\). From \cite[Proposition 1.36]{hat}, it follows that \te a two sheeted covering \(\pi:\tx \ra X\) \st \(\tx\) is connected and if 
	\[\pi_{*}:\pi_1(\tx)\ra \pi_1(X),\]
	\(\text{im}(\pi_{*})=\ker(\Ph_*)\). By \cite{hat}*{Proposition 1.33}, \(\Ph\) has a lift 
	\[\tp: \tx \ra (\cm, \cF) \quad \text{such that} \quad \del \circ \tp = \Ph \circ \pi.\]
	Moreover, \(\tp\) is \(\bbz_2\)-equivariant, i.e. if \(T:\tx \ra \tx\) is the deck transformation, then for all \(x\in \tx\),
	\[\tp(T(x))=M-\tp(x).\]	
	
	Let us recall the notation from Section \ref{sec.2.4}.
	\[\stx=\text{ suspension of } \tx= \frac{\tilde{X} \times [-1,1]}{\sim}.\]
	By the abuse of notation, an element of $S\tx$ will be denoted by a pair $(x,t)$ with $x \in \tx, t \in [-1,1]$.
	\[S\tx = C_{+}\cup C_{-}, \quad \text{where} \quad C_{+} = \frac{\tx \times [0,1]}{\sim} \; , \; C_{-} = \frac{\tx \times [-1,0]}{\sim}.\]
	Since $\tx$ is connected, $S\tx$ is simply-connected. There is a free $\mathbb{Z}_2$ action on $S\tx$ given by $(x,t) \mapsto (T(x),-t)$. Let $Y$ denote the quotient of $S\tx$ with respect to this $\mathbb{Z}_2$ action and $\rh:S\tx \ra Y$ be the covering map. \(\tx\) naturally sits inside \(\stx\) by the map \(x \mapsto (x,0)\); hence, $X$ also naturally sits inside $Y$. We will identify \(\tx\) with its image under \(\tx \hookrightarrow S\tx\) and \(X\) with its image under \(X \hookrightarrow Y\). 
	
	To illustrate the above constructions, let us consider the following example. If $X=\mathbb{RP}^k$, then $\tx = S^k,\; S\tx = S^{k+1}$, $C_+$ and $C_{-}$ are respectively the upper hemisphere and the lower hemisphere of \(S^{k+1}\), $Y=\mathbb{RP}^{k+1}$.
	
\vspace{1\baselineskip}
\textbf{Part 2.} We begin with the following proposition from \cite{afp}. For the shake of completeness, we also include its proof (following \cite{afp}).
\begin{pro}[\cite{afp}*{Proposition 3.38}]
	\label{pro intersection of cac sets}
For an \(\cH^{n+1}\)-measurable set \(E\subset M\) and an open set \(U \subset M\), let us use the notation
\[P(E,U)=\sup \left\{\int_{E \cap U}\operatorname{div } \om\; d\mathcal{H}^{n+1}: \om \in \mathfrak{X}^1_c(U) \text{ and } \|\om\|_{\infty}\leq 1 \right\}.\] 
Here \(\mathfrak{X}^1_c(U)\) is the space of compactly supported \(C^1\) vector-fields on \(U\). Suppose \(E,F\in \cm\). Then
\[P(E \cap F,U) \leq P(E,U) + P(F,U).\]
Hence, in particular, \(E \cap F \in \cm\).
\end{pro}
\begin{proof}
By \cite{MPPP}*{Proposition 1.4}, there exist sequences \(\{f_j\}_{j=1}^{\infty},\{g_j\}_{j=1}^{\infty}\subset C^{\infty}(U)\) \st \(0\leq f_j,g_j\leq 1\) for all \(j \in \bbn\),
\begin{align*}
&f_j \ra \ch_E|_U \; \text{ and } \; g_j \ra \ch_F|_U \; \text{ in }\; L^1(U);\\
&P(E,U)=\lim_{j \ra \infty} \int_{U}\|Df_j\|\; d\cH^{n+1} \; \text{ and } \; P(F,U)=\lim_{j \ra \infty} \int_{U}\|Dg_j\|\; d\cH^{n+1}.
\end{align*}
By the dominated convergence theorem, 
\[f_jg_j \ra \ch_E\ch_F|_U=\ch_{E \cap F}|_U \; \text{ in } \;L^1(U).\]
\tf
\begin{align*}
	P(E \cap F , U)&\leq \liminf_{j \ra \infty}\int_{U}\|D(f_jg_j)\|\; d\cH^{n+1}\\&\leq \liminf_{j \ra \infty}\int_{U}\|Df_j\|\; d\cH^{n+1}+ \liminf_{j \ra \infty}\int_{U}\|Dg_j\|\; d\cH^{n+1}\\&=P(E,U)+P(F,U).
\end{align*}
\end{proof}
The map $\tp:\tx \ra (\cm,\cF)$ (defined in Part 1) can be extended to a map $\ts':S\tx \ra (\cm,\cF)$ as follows. Let us define
	\[\tz =\{(x,t)\in \stx : |t|\leq 1/2\}\subset \stx.\]
    Let $\Lambda:[0,1] \ra (\cm,\cF)$ be \st $\Lambda(0)=M$, $\Lambda(1)=\emptyset$, \(\del \circ \La :[0,1] \ra \znf\) has no concentration of mass and 
	\begin{equation}\label{E:mass of del Lambda}
		\sup_{t \in [0,1]}\mathbf{M}(\partial \Lambda(t))\leq W_0+\et/2-\delta.
	\end{equation}
	 We define	
	\begin{equation}\label{E:defn of tilde Psi}
	\ts'(x,t)=
	\begin{cases}
	\tp(x) & \text{ if } (x,t)\in \tz; \\
	\tp(x) \cap \Lambda(2t-1) & \text{ if } (x,t) \in C_{+}\setminus\tz; \\
	M- \big(\tp(T(x)) \cap \Lambda(-1-2t)\big) & \text{ if } (x,t) \in C_{-}\setminus\tz.
	\end{cases}
	\end{equation}
By Proposition \ref{pro intersection of cac sets}, \(\ts'(x,t)\in \cm\) for all \((x,t)\in S\tx.\)
\begin{clm}
	$\ts':S\tx \ra \cm$ is continuous in the flat topology.
\end{clm}
\begin{proof}
	It is enough to prove that \((x,t) \mapsto \tp(x) \cap \Lambda(2t-1)\) is \cts in \(\cF\) for \((x,t)\in C_{+}\setminus\tz\). We note that for arbitrary sets \(A_1\), \(A_2\), \(B_1\), \(B_2\),
	\[\left(A_1\cap A_2\right)\De\left(B_1\cap B_2\right) \subset \left(A_1 \; \De \; B_1\right)\cup \left(A_2 \;\De \;B_2\right),\]
	where \(A\;\De\; B= (A \setminus B)\cup (B \setminus A)\) denotes the symmetric difference of \(A\) and \(B\). Hence, if \((x_1,t_1),(x_2,t_2) \in C_{+}\setminus\tz\),
\[\big(\tp(x_1) \cap \Lambda(2t_1-1)\big)\De \big(\tp(x_2) \cap \Lambda(2t_2-1)\big) \subset \big(\tp(x_1)\;\De\;\tp(x_2)\big) \cup \big(\Lambda(2t_1-1) \; \De \; \Lambda(2t_2-1)\big),\]
which implies
\begin{align}
\label{e pf of the claim}
 \cF\Big(\big(\tp(x_1) \cap \Lambda(2t_1-1)\big)& - \big(\tp(x_2) \cap \Lambda(2t_2-1)\big)\Big)\nonumber\\
& \leq \cF\big(\tp(x_1)-\tp(x_2)\big) + \cF\left(\Lambda(2t_1-1) -\Lambda(2t_2-1)\right).
\end{align}
Since \(x \mapsto \tp(x)\) and \(t \mapsto \La(t)\) are continuous in the flat topology, by \eqref{e pf of the claim}, \((x,t) \mapsto \tp(x) \cap \Lambda(2t-1)\) is also continuous in the flat topology.
\end{proof}
	As $\partial\ts' (x,t) = \partial\ts' (T(x), -t)$, $\ts'$ descends to a \cts map $\Psi' : Y \ra \znf$, i.e. \(\Ps'\circ \rh = \del\circ\tilde{\Psi}'\). Furthermore, $\Ps'|_X=\Ph.$	
\begin{clm}\label{cl ncm}
$\Psi' : Y \ra \znf$ has no concentration of mass.	
\end{clm}
\begin{proof}
	If \(E \in \cm\) and \(U \subset M\) is open, \(P(E,U)=\|\del E\|(U).\) By Proposition \ref{pro intersection of cac sets}, for \(r>0\) and \(p\in M\),
\begin{equation*}
\sup_{y \in Y}\|\Ps'(y)\|(B(p,r)) \leq \sup_{x \in X} \|\Ph(x)\|(B(p,r)) + \sup_{t \in [0,1]} \|\del \La(t)\|(B(p,r)).
\end{equation*}
This finishes the proof of the claim as the maps \(\Ph\) and \(\del \circ \La\) have no concentration of mass.
\end{proof}	
By the interpolation theorems of Marques and Neves \cite{MN_Willmore}*{Theorem 13.1, 14.1}, \cite{mn_ricci_positive}*{Theorem 3.9, 3.10}, the above Claim \ref{cl ncm} implies that there exists a \cts map \(\Ps:Y \ra \zn\) \st $\Ps$ is homotopic to \(\Ps'\) in the \(\cF\)-topology and for all \(y \in Y\),
\begin{equation}
	\label{e. mass diff of Ps and Ps'}
	\bM(\Ps(y))\leq \bM(\Ps'(y))+\de/2.
\end{equation}

\vspace{1\baselineskip}
\textbf{Part 3.} We will prove that \(\Ps:Y\ra\zn\), constructed above, is a \((k+1)\)-sweepout.

Let \(\overline{\la}\) be the generator of \(H^*(\znf,\bbz_2)\) as in Section \ref{sec.2.3} and \(\la=\Ps^*\overline{\la}\). We need to show that \(\la^{k+1}\neq 0 \in H^{k+1}(Y,\bbz_2)\). Since \(\stx\) is simply-connected and $Y$ is the quotient of \(\stx\) under the \(\mathbb{Z}_2\) action, \(\pi_1(Y)=\bbz_2\). Hence, \(H_1(Y,\bbz_2)=H^1(Y,\bbz_2)=\bbz_2\) as well.

Let \(\iota:X \ra Y\) be the inclusion map. Since \(\Ps'|_X=\Ph\), \(\Ps\) is homotopic to \(\Ps'\) in the \(\cF\)-topology and \(\Ph\) is a \(k\)-sweepout, \(\iota^*\la=\Ph^*\overline{\la}\neq 0 \in H^1(X,\bbz_2)\). Therefore, \(\la\) is the unique non-zero element of \(H^1(Y,\bbz_2)\). Hence, if \(\ga:S^1 \ra Y\) is a non-contractible loop, \(\la.[\ga]=1\).

To prove \(\la^{k+1}\neq 0 \in H^{k+1}(Y,\bbz_2)\), it is enough to find \(\al \in H_{k+1}(Y,\bbz_2)\) \st \(\la^{k+1}.\al=1\). By \cite{hat}*{Proof of Proposition 2B.6}, \te a map (called the \textit{transfer homomorphism}) \(\tau_*:H_k(X,\bbz_2)\ra H_k(\tx,\bbz_2)\), which is induced by the chain map \(\ta :C_k(X,\bbz_2)\ra C_k(\tx,\bbz_2)\), defined as follows. If \(u:\De^k\ra X\), \(u\) has precisely two lifts \(\tilde{u}_1,\tilde{u}_2:\De^k \ra \tx\). We define \(\tau(u)=\tilde{u}_1+\tilde{u}_2\). 

As \(\Ph\) is a \(k\)-sweepout, \((\iota^*\la)^k=(\Ph^*\overline{\la})^k\neq 0 \in H^k(X,\bbz_2)\). \tf \te \(\si \in H_k(X,\bbz_2)\) \st \((\iota^*\la)^k.\si=1\). Let \(\si = [\sum_{i=1}^{I}\si_i]\), where each \(\si_i:\De^k\ra X\) is a singular \(k\)-chain. Let \(\ta(\si_i)=\tsia+\tsib\), where \(\tsia,\tsib :\De^k \ra \tx\) are lifts of \(\si_i\). Since
\(\del\big(\sum_i \si_i \big)=0 \)
and \(\ta\) is a chain map, we have
\begin{equation}\label{E: tau sigma is cycle}
\del\Big(\sum_i \tsia+\tsib \Big)=0
\end{equation}
as well. \(C\De^k\) can be naturally identified with \(\De^{k+1}\) such that the collapsed image of \(\De^k\times\{1\}\) is the \((k+1)\)-st vertex \(v_{k+1}\) (see Section \ref{sec.2.4} for the notation). Let us also recall from Part 1 that \(C_+\) is a cone over \(\tx\) ; let \(y_0 \in C_+\) be the collapsed image of \(\tx \times \{1\}\). For \(s=1,\dots,2I\), we consider 
\[C\tilde{\si}_s:C\De^k(=\De^{k+1})\ra C\tx(=C_+); \; \; C\tilde{\si}_s(v_{k+1})=y_0.\]
Using equation \eqref{E: tau sigma is cycle}, one can check that 
\begin{equation}\label{c tau sigma is cycle}
	\del\Big(\sum_i C\tsia + C\tsib \Big) = \Big(\sum_i \tsia+\tsib \Big).
\end{equation}
As \(\tsia,\tsib\) are lifts of \(\si_i\), \(\rh_\#\tsia=\rh_\#\tsib=\si_i\). Hence,
\begin{equation}\label{E: alpha is cycle}
\del\Big(\sum_i \rh_\#C\tsia + \rh_\#C\tsib \Big) = \rh_\#\Big(\sum_i \tsia+\tsib \Big)=0 \in C_k(X,\bbz_2).
\end{equation}
Let us define
\begin{equation}
\label{E: defn of alpha}
\al=\Big[\sum_i \rh_\#C\tsia + \rh_\#C\tsib \Big]\in H_{k+1}(Y,\bbz_2).
\end{equation}
We claim that \(\la^{k+1}.\al=1\). Let \(\la=[\ell], \; \ell \in C^{k+1}(Y,\bbz_2),\; \text{i.e. }\ell:C_{k+1}(Y,\bbz_2)\ra\bbz_2\). By the definition of the cup product,
\begin{align}
&\ell^{k+1}\Big(\rh_\#C\tsia + \rh_\#C\tsib \Big) \nonumber \\
&=\ell^k\Big(\rh_\#C\tsia\big|_{\De^k}\Big)\ell\Big(\rh_\#C\tsia\big|_{[v_k,v_{k+1}]}\Big) + \ell^k\Big(\rh_\#C\tsib\big|_{\De^k}\Big)\ell\Big(\rh_\#C\tsib\big|_{[v_k,v_{k+1}]}\Big) \nonumber \\
&=\ell^k(\si_i)\ell\Big(\rh_\#C\tsia\big|_{[v_k,v_{k+1}]}+\rh_\#C\tsib\big|_{[v_k,v_{k+1}]}\Big). \label{e-l.loop 1}
\end{align}
In the second equality we have used the fact that \[C\tsia|_{\De^k}=\tsia \; \text{ and } \; C\tsib|_{\De^k}=\tsib.\] 
We note that \[C\tsia|_{[v_k,v_{k+1}]} \text{ is a curve joining } \tsia(v_k) \text{ and } \tsia(v_{k+1})=y_0,\]
and  
\[C\tsib|_{[v_k,v_{k+1}]} \text{ is a curve joining } \tsib(v_k) \text{ and } \tsib(v_{k+1})=y_0.\]

 As \(\tsia,\tsib\) are lifts of \(\si_i\), \(T\left(\tsia(v_k)\right)=\tsib(v_k)\). Hence, \[\left(\rh_\#C\tsia\big|_{[v_k,v_{k+1}]}+\rh_\#C\tsib\big|_{[v_k,v_{k+1}]}\right)\] is a non-contractible loop in \(Y\). \tf 
\begin{equation}
\ell\left(\rh_\#C\tsia\big|_{[v_k,v_{k+1}]}+\rh_\#C\tsib\big|_{[v_k,v_{k+1}]}\right)=1 \label{e-l.loop 2}
\end{equation}
Using \eqref{e-l.loop 1} and \eqref{e-l.loop 2}, we conclude that
\begin{equation*}
\la^{k+1}.\al=\sum_i \ell^{k+1}\Big(\rh_\#C\tsia + \rh_\#C\tsib \Big) \nonumber  = \sum_i \ell^k(\si_i)= (\iota^*\la)^k.\si = 1.
\end{equation*}
\medskip

\textbf{Part 4.} Since \(\Ps\) is a \((k+1)\)-sweepout, using the argument as in Part 1, \(\Ps\) has a lift 
\[\ts: S\tx \ra (\cm, \cF) \quad \text{such that} \quad \del \circ \ts = \Ps \circ \rh.\]
\(\ts\) is \(\bbz_2\)-equivariant, i.e.
\[ \ts(T(x),-t)=M-\ts(x,t),\]
for all \((x,t)\in S\tx\). Moreover, by the definition of the \(\bF\)-\mt on \(\cm\) and \(\mathcal{Z}_n(M;\bbz_2)\), continuity of \(\Ps:Y \ra \zn\) implies that \(\ts:S\tx\ra (\cm,\bF)\) is also continuous. Let $\tilde{\Pi}$ be the \((\stx,\tz)\)-homotopy class of \(\ts\). We will show that 
\begin{equation}
\label{e-Pi is non-trivial}
\bL^c(\tilde{\Pi})>\sup_{y \in \tz }\ac(\ts(y)).
\end{equation}

The proof closely follows \cite{zhou2}*{Proof of Lemma 5.8}. By \eqref{e. mass diff of Ps and Ps'}, \eqref{E:defn of tilde Psi} and \eqref{e. mass of Phi in X},
\begin{align}
\sup_{y \in \tz }\ac(\ts(y))\leq \sup_{y \in \tz }\bM(\del\ts(y))&\leq\sup_{y \in \tz }\bM(\del\ts'(y))+\de/2 \nonumber \\
&=\sup_{x \in \tx}\bM(\del \tp(x))+\de/2 \leq \om_k+\de.\label{e. estimating ac on tz}
\end{align}
Let \(\{\ts_i:\stx \ra (\cm,\bF)\}_{i=1}^{\infty}\) be an arbitrary element in \(\tilde{\Pi}\). There exists \(H_i:\stx \times [0,1] \ra (\cm,\bF)\) \st \(H_i(-,0)=\ts\), \(H_i(-,1)=\ts_i\) and
\begin{equation}
\limsup_{i \rightarrow \infty} \sup \left\lbrace\bF(H_i(y,s),\ts(y)): y \in \tz, s\in [0,1]\right\rbrace =0. \label{e. deviation goes to 0}
\end{equation}
We define \(\{\ts_i^*:C_+ \ra (\cm,\bF)\}_{i=1}^{\infty}\) as follows.
\begin{equation}
\ts^*_i(y)=
\begin{cases}
\ts(y) & \text{ if } y \in \tx; \\
H_i(y,2t) & \text{ if } y=(x,t) \in \tz \cap C_+ ;\\
\ts_i(y) & \text{ if } y \in C_+ \setminus \tz.
\end{cases}
\label{e. defn of ts_i*}
\end{equation}
\(\rh|_{C_+}:C_+ \ra Y\) is a quotient map. Since \(\ts\) is \(\bbz_2\)-equivariant, \(\ts^*_i\) descends to a map \(\Ps^*_i : Y \ra \zn\) \st \(\Ps^*_i \circ (\rh |_{C_+})=\del \circ \ts^*_i \). Moreover, $\Ps_i^*$ is homotopic to $\Ps$ (in the \(\bF\)-topology) and hence, is a \((k+1)\)-sweepout. \(\Ps^*_i\) has no concentration of mass as it is continuous in the \(\bF\)-topology. \tf 
\begin{equation}
\label{e. mass of ts_i*}
\sup_{y \in C_+}\bM(\del\ts^*_i(y))=\sup_{y \in Y}\bM(\Ps^*_i(y))\geq \om_{k+1}.
\end{equation}
However, as noted in \eqref{e. estimating ac on tz},
\begin{equation*}
\sup_{y \in \tz }\bM(\del\ts(y))\leq \om_k+\delta < \om_{k+1}\; (\text{using } \eqref{E:defn of delta});
\end{equation*}
hence, \eqref{e. deviation goes to 0} and \eqref{e. defn of ts_i*} imply that if \(i\) is sufficiently large,
\begin{equation}
\label{e. mass of ts^*_i in tz}
\sup_{y \in \tz\cap C_{+}}\bM(\del\ts^*_i(y)) < \om_{k+1}
\end{equation}
as well. Since \(\ts^*_i\) coincides with \(\ts_i\) on \(C_+\setminus\tz\), \eqref{e. mass of ts_i*} and \eqref{e. mass of ts^*_i in tz} imply that for \(i\) sufficiently large,
\begin{align}
\sup_{y \in C_+ \setminus \tz }\bM(\del\ts_i(y))\geq \om_{k+1} & \Longrightarrow \sup_{y \in \stx }\bM(\del\ts_i(y))\geq \om_{k+1}\nonumber \\
& \Longrightarrow \sup_{y \in \stx }\ac(\ts_i(y))\geq \om_{k+1}-c\vol(g).\label{e. ac ts i}
\end{align}
Since \eqref{e. ac ts i} holds for any \(\{\ts_i\}_{i=1}^{\infty}\in \tilde{\Pi}\), using \eqref{E:defn of delta} and \eqref{e. estimating ac on tz} we conclude that
\begin{equation}
\bL^c(\tilde{\Pi})\geq \om_{k+1}-c\vol(g)>\om_k+\delta\geq \sup_{y\in\tz}\ac(\ts(y)).\label{e. lower bound of L^c Pi}
\end{equation}

\textbf{Part 5.} By Theorem \ref{min-max thm} and \eqref{e-Pi is non-trivial}, there exists \(\Om \in \cm\) \st $\del \Om$ is a closed $c$-CMC (\wrt the inward unit normal) \hy with optimal regularity and $\ac(\Om)=\bL^c(\tilde{\Pi})$. \eqref{e. lower bound of L^c Pi} implies that
\begin{equation}
\ac(\Om)=\bL^c(\tilde{\Pi}) >\om_k+\de>\om_k.
\end{equation}
Moreover, 
\begin{align*}
\ac(\Om)=\bL^c(\tilde{\Pi}) &\leq \sup_{(x,t) \in \stx }\ac(\ts(x,t))\\
& \leq \sup_{(x,t) \in \stx }\bM(\del\ts(x,t))\\
& \leq \sup_{(x,t) \in \stx }\bM(\del\ts'(x,t))+\de/2\;\;(\text{by } \eqref{e. mass diff of Ps and Ps'})\\
& \leq \sup_{x\in\tx}\bM(\del\tp(x))+\sup_{t\in[0,1]}\bM(\del\La(t))+\de/2\;\;(\text{by \eqref{E:defn of tilde Psi} and Proposition }\ref{pro intersection of cac sets})\\
& \leq \om_k+\de/2+W_0+\et/2-\de+\de/2 \;\;(\text{by }\eqref{e. mass of Phi in X} \text{ and } \eqref{E:mass of del Lambda}) \\
& <\om_k+W_0+\eta.
\end{align*}
This finishes the proof of Theorem \ref{main thm}.
\end{proof}
\vspace{1\baselineskip}
\begin{proof}[Proof of Theorem \ref{main cor}] We start with the following lemma.
\begin{lem}
	\label{ac omega1 neq ac omega2} Let \(\Om_1,\;\Om_2 \in \cm\) and \(c>0\) \st for \(i=1,2\), \( \del \Om_i \) is a closed $c$-CMC (\wrt the inward unit normal) \hy with optimal regularity and \(\ac(\Om_1)\neq \ac(\Om_2)\). Then \(\del\Om_1 \neq \del \Om_2\).
\end{lem}
\begin{proof}
	Suppose \(\Om_1=\Om_2=\Si\). Since \(\ac(\Om_1) \neq \ac(\Om_2)\), \(\Om_1 \neq \Om_2\). Hence, by the constancy theorem, \(\Om_2 = M \setminus \Om_1\); this is impossible as by our hypothesis, the mean curvature vector of \(\Si\) points inside both \(\Om_1\) and \(\Om_2\).
\end{proof}

Since \(\{\om_p\}_{p=1}^{\infty}\) is a non-decreasing sequence of positive real numbers converging to \(+\infty\), for every \(m \in \bbn\), \te a subsequence \(\{\om_{p_i}\}_{i=1}^{m}\) \st 
\[\om_{p_i}<\om_{p_i+1} \quad \text{ and }\quad \om_{p_i}+W_0 < \om_{p_{i+1}},\]
for all \(i = 1,\dots,m.\) Let us define
\begin{equation}
\label{e. defn of c*(m)}
c^{*}(m) = \frac{1}{\vol(g)}\min\left\{\om_{p_i+1}-\om_{p_i}:i=1,\dots,m\right\}.
\end{equation}
If \(0<c<c^{*}(m)\), applying Theorem \ref{main thm} for \(k=p_i\), \(i=1,\dots,m\), we conclude that there exists \(\Om_i \in \cm\), \st \( \del \Om_i \) is a closed $c$-CMC (\wrt the inward unit normal) \hy with optimal regularity and 
\[\om_{p_i}<\ac(\Om_i)<\om_{p_{i+1}}.\]
\tf using Lemma \ref{ac omega1 neq ac omega2}, \(\{\del \Om_i\}_{i=1}^m\) are distinct \(c\)-CMC hypersurfaces. 

To prove the more quantitative statement, we will use the growth estimate of the volume spectrum \eqref{e. sublinear growth}. From the proof of Theorem \ref{main thm}, one can obtain the following inequality.
\begin{equation}
	\om_{k+1}\leq \om_{k}+W_0 \quad \forall k\in \bbn. \label{e.subadditive property of omega_k}
\end{equation}
Let us briefly repeat the argument here. We fix \(\de>0\). There exists a \(k\)-sweepout \(\Ph:X \ra \znf\) with no concentration of mass \st \(\sup\left\{\bM(\Ph(x)):x \in X\right\}\leq \om_k+\de\). Using the notation from the Proof of Theorem \ref{main thm}, Part 1, \te a double cover \(\pi:\tilde{X}\ra X\) \st \(\Ph\) has a lift \(\tp:\tx \ra (\cm,\cF)\). Let \(\La:[0,1]\ra (\cm,\cF)\) be \st \(\La(0)=M\), \(\La(1)=\emptyset\), \(\del \circ \La :[0,1] \ra \znf\) has no concentration of mass and \(\sup\left\{\bM(\del\La(t)):t \in [0,1]\right\} \leq W_0+\de\). We define the map \(\ts:\stx \ra (\cm,\cF)\) as follows (cf. \eqref{E:defn of tilde Psi}).
\begin{equation*}
\ts(x,t)=
\begin{cases}
\tp(x) \cap \Lambda(t) & \text{ if } (x,t) \in C_{+}; \\
M- \big(\tp(T(x)) \cap \Lambda(-t)\big) & \text{ if } (x,t) \in C_{-}.
\end{cases}
\end{equation*}
\(\ts\) descends to a map \(\Ps:Y\ra \znf\). As proved in Part 2 (using Proposition \ref{pro intersection of cac sets}), \(\Ps\) has no concentration of mass. Moreover, as argued in Part 3, \(\Ps\) is a \((k+1)\)-sweepout; hence, 
\begin{equation}\label{sa}
\sup\left\{\bM(\Ps(y)):y \in Y\right\}\geq \om_{k+1}.
\end{equation}
Further, as estimated in Part 5 (using Proposition \ref{pro intersection of cac sets}),
\begin{equation}\label{sb}
\sup\left\{\bM(\Ps(y)):y \in Y\right\}\leq \om_k+W_0+2\de.
\end{equation}
Since \(\de > 0\) is arbitrary, \eqref{sa} and \eqref{sb} imply \eqref{e.subadditive property of omega_k}.

As a consequence of \eqref{e.subadditive property of omega_k}, for all \(N\geq \om_1\), there exist \(q,r\in \bbn\) \st \(\om_q \in [N,N+W_0) \) and \(\om_r \in [N+2W_0,N+3W_0)\). Using \eqref{e. sublinear growth},
\[K_0r^{\frac{1}{n+1}}\leq \om_r< N+3W_0\quad \Longrightarrow r < \left(\frac{N+3W_0}{K_0}\right)^{n+1}. \]
Hence, there exists \(s \in \bbn\), \(q\leq s \leq r\), \st 
\begin{equation*}
\om_{s+1}-\om_s\geq \frac{\om_r-\om_q}{r-q}> W_0\left(\frac{N+3W_0}{K_0}\right)^{-(n+1)}.
\end{equation*}
\tf by Theorem \ref{main thm}, if
\[0<c\vol(g)\leq W_0\left(\frac{N+3W_0}{K_0}\right)^{-(n+1)},\]
there exists \(\Om \in \cm\), \st \(\del \Om\) is a closed $c$-CMC (\wrt the inward unit normal) \hy with optimal regularity and
\[N\leq\om_s < \ac(\Om)<\om_s+W_0+(N+3W_0-\om_s)= N+4W_0.\]
We note that \(W_0\geq \om_1\). Setting \(N=(4i-3)W_0\), \(i=1,\dots ,m\) in the above argument, we conclude that if 
\[0<c\vol(g)\leq W_0\left(\frac{4iW_0}{K_0}\right)^{-(n+1)},\]
\te \(\Om_i \in \cm\), \st \(\del \Om_i\) is a closed $c$-CMC (\wrt the inward unit normal) \hy with optimal regularity and
\[(4i-3)W_0 < \ac(\Om_i)<(4i+1)W_0.\]
Further, from the one parameter min-max construction of Zhou and Zhu \cite{zz1}, for any \(c>0\) \te \(\Om_0 \in \cm\) \st \(\del \Om_0\) is a closed $c$-CMC (\wrt the inward unit normal) \hy with optimal regularity and
\(\ac(\Om_0)\leq W_0\).
\tf if
\begin{equation} \label{asymtotic formula}
0<c\vol(g)\leq W_0\left(\frac{4mW_0}{K_0}\right)^{-(n+1)},
\end{equation}
by Lemma \ref{ac omega1 neq ac omega2}, there exist at least \((m+1)\) many closed \(c\)-CMC \hys in \((M,g)\). This implies the number of closed \(c\)-CMC \hys in \((M,g)\) is at least
\[\left\lfloor\frac{K_0}{4W_0}\left(\frac{W_0}{c\vol(g)}\right)^{\frac{1}{n+1}}\right\rfloor+1, \]
where for \(x \in \bbr^+\), \(\lfloor x \rfloor\) is the largest (non-negative) integer \(\leq x\) and \(K_0\) is the constant appearing in \eqref{e. sublinear growth}.
\end{proof}

\begin{rmk}
For \(m \in \bbn\), the value of \(c^*(m)\), obtained from \eqref{e. defn of c*(m)} for an optimal subsequence \(\{\om_{p_i}\}_{i=1}^{m}\), will be greater than or equal to the value, obtained from the asymptotic formula \eqref{asymtotic formula}.
\end{rmk}

\begin{ex}
	Let us consider \(S^3\) equipped with the round metric \(g_0\). By the one parameter min-max construction of Zhou and Zhu \cite{zz1}, for any \(c \in \bbr^+\) \te \(\Om \in \mathcal{C}(S^3)\), \st \(\del \Om \) is a closed $c$-CMC (\wrt the inward unit normal) \hy and \(\ac(\Om)\leq W_0(S^3, g_0)=4\pi\). As proved by Nurser \cite{nurser}, for the \mt \(g_0\), \(\om_1=\dots=\om_4=4\pi\) and \(\om_5 = 2\pi^2\). \tf by Theorem \ref{main thm}, if
	\[0<c<\frac{\om_5 - \om_4}{\vol(g_0)}= \frac{2\pi^2 - 4\pi}{2\pi^2} = 1 - \frac{2}{\pi} \approx \frac{1}{3},\]
	\te \(\Om'\in \mathcal{C}(S^3)\) \st \(\del \Om' \) is a closed $c$-CMC (\wrt the inward unit normal) \hy and \(\ac(\Om')>\om_4(S^3, g_0)= 4\pi\). Thus, for the \mt \(g_0\) (with the notation as in Theorem \ref{main cor}), \(c^{*}(2)\) can be taken to be \(1-2/\pi\). Moreover, the widths are continuous functions of the \mt \cite{imn,mns} and \(W_0\) coincides with \(\om_1\) when the ambient metric has positive Ricci curvature \cite{zhou_ricci_positive,zhou.higher.dim}. \tf if \(\tilde{g}\) is a small perturbation of \(g_0\), \(c^{*}(2)\) for \(\tilde{g}\) can be taken to be close to \(1-2/\pi\).
\end{ex}

\bibliographystyle{amsalpha}
\bibliography{multiplecmc}

\end{document}